\documentclass[letterpaper,10pt,conference]{ieeeconf}
\IEEEoverridecommandlockouts                         
\usepackage{amsfonts}

\usepackage{amsmath}
\usepackage{graphicx}
\usepackage{cite}
\usepackage{color}
\usepackage{verbatim}

\usepackage{hyperref}
\hypersetup{
    colorlinks=false,
    pdfborder={0 0 0},
}

\title{The Quadratic-Quadratic Regulator Problem:\\ Approximating feedback controls for quadratic-in-state nonlinear systems}
\author{Jeff Borggaard and Lizette Zietsman%
\thanks{This research was partially supported by the National Science Foundation under contract DMS-1819110 and the authors gratefully acknowledge the support of the Institute for Mathematics and its Applications (IMA), where much of the work was carried out during the annual program on Control Theory and its Applications.}%
\thanks{Jeff Borggaard and Lizette Zietsman are
with the Interdisciplinary Center for Applied Mathematics and the Department of Mathematics at
Virginia Tech, Blacksburg, VA 24061, USA.\newline
{\tt\footnotesize jborggaard@vt.edu, lzietsma@vt.edu}.}
}

\begin{document}
\maketitle

\begin{abstract}
Feedback control problems involving autonomous quadratic systems are prevalent, yet there are only a limited number of software tools available for approximating their solution due to the complexity of the problem.  This paper represents a step forward in the special case where both the state equation and the control costs are quadratic.  As it represents the natural extension of the {\em linear-quadratic regulator} (LQR) problem, we describe this setting as the {\em quadratic-quadratic regulator} (QQR) problem.  This is significantly more challenging and holds the LQR as special case that must be solved along the way.  We describe an algorithm that exploits the structure of the QQR problem that arises when implementing Al'Brekht's method. This approach is amenable to feedback laws with low degree polynomials but have a relatively modest model dimension that could be achieved by modern model reduction methods.  This problem has an elegant formulation and a solution that introduces several linear systems where the structure suggests modern tensor-based linear solvers.  We demonstrate this algorithm on a suite of random test problems then apply it to a distributed parameter control problem that fits the QQR framework.  Comparisons to linear feedback control laws show a modest benefit using the QQR formulation.
\end{abstract}

\section{Motivation}
Linear feedback control of autonomous nonlinear systems, such as those describing the behavior of fluids, can be sufficient to achieve stabilization--even for an unstable steady-state solution~\cite{bergmann2005orc,benner2016rsl,borggaard2010lfc,borggaard2014mrd}.  
There is a shortage of software tools for nonlinear problems in control and systems theory.  The general Matlab Nonlinear Systems Toolbox (NST) by Krener~\cite{krener2015NST} takes a broad step toward delivering useful tools for a number of important problems.  Since we inherently encounter the {\em curse of dimensionality} in these problems, there is also a need to develop specialized tools for important classes of problems.  This paper addresses this by specifically solving the quadratic-quadratic regulator problem: minimizing a quadratic cost subject to a state equation with a quadratic nonlinearity. 

For example, linear feedback laws found by solving the linear-quadratic regulator (LQR) problem compute the linear feedback law as the solution to a single algebraic Riccati equation and have the property that the linear portion of the nonlinear system becomes stable~\cite{barbu2007les,barbu2001fip,barbu2006tbs,raymond2006fbs}.  Unfortunately, for nonlinear systems, this only guarantees local stability.  Thus the ability of linear feedback to stabilize the steady-state solution depends on the initial condition, which must be sufficiently close to the steady-state.  An alternative would be to develop nonlinear feedback control laws that could offer the ability to  expand the radius of convergence (shown with a simple example in \cite{borggaard2018computation}).  However, these require us to approximate solutions to the Hamilton-Jacobi-Bellman (HJB) equations, e.g.~\cite{kunisch2004hjb,breiten2019feedback}.  The HJB equations are notoriously complex in the general case.  Nevertheless, if one considers the quadratic-quadratic regulator (QQR) problem, with autonomous quadratic state equations and a quadratic control objective, there is sufficient structure in polynomial approximations based on Al'Brekht's method \cite{navasca2000solution} for polynomial feedback laws to be computable for modest problem sizes.  The QQR problem also happens to be exactly what is needed to solve discretized versions of distributed parameter control problems where the nonlinearity is quadratic (such as the Navier-Stokes equations used as our motivation above).  This is particularly true when linear feedback laws are being based on LQR problems. As in the LQR case, suitable model reduction methods \cite{aubry1988dcs,holmes1996tcs,ahmad2015nim} are essential to forming a solution methodology for distributed parameter control problems with quadratic nonlinearities. First of all, the Riccati equation is still needed to compute the linear term \cite{singler2008alr,benner2013esl,singler2016pod} and the curse-of-dimensionality still appears with higher-order polynomial approximations of the feedback law.  

In this paper, we briefly outline the HJB equations, the QQR problem, and polynomial approximations to the value function and the feedback control operators.  Our formulation leads to a sequence of linear systems in Kronecker product form after an initial solution to the algebraic Riccati equation.  As we shall see, a na\"ive construction of these matrices and other terms would quickly become prohibitive.  However, the structure lends itself to newly developed recursive tensor linear algebra that avoids assembly and other taxing of computer memory.  We present a numerical study with a set of randomly selected control problems to compare solutions obtained by Krener's NST~\cite{krener2015NST}, direct assembly and solution to the Kronecker system, and the recursive tensor-based algorithm using the tensor toolbox~\cite{kolda2009tensor} and a recursive blocked algorithm for systems with a special Kronecker sum form~\cite{chen2019recursive}.

\section{Background}
For simplicity of exposition, we describe the nonlinear optimal control problem and its computational challenges for systems modeled by autonomous systems of ordinary differential equations.  This can be widely found in the literature and we are reintroducing it here to set up our notation.  The problem is to find a control ${\bf u}(\cdot)\in L_2(0,\infty;\mathbb{R}^{m})$ that solves
\begin{equation}
\label{eq:oc}
  \min_{\bf u} J({\bf x},{\bf u}) = \int_0^\infty \ell({\bf x}(t),{\bf u}(t))\ dt,
\end{equation}
where $\ell : \mathbb{R}^n \times \mathbb{R}^{m} \longrightarrow [0,\infty)$ is a prescribed control objective, and minimization occurs subject to
\begin{equation}  \label{eq:full}
  \dot{{\bf x}}(t) = {\bf A}{\bf x}(t) + {\bf B} {\bf u}(t) + {\bf f}({\bf x}(t),{\bf u}(t)),
\end{equation}
from ${\bf x}(0) = {\bf x}_0\in\mathbb{R}^n$, where 
${\bf A} \in\mathbb{R}^{n\times n} $ and ${\bf B} \in \mathbb{R}^{n\times m}$ are constant matrices and $ \mathbf{f}: \mathbb{R}^{n}\times\mathbb{R}^{m}\longrightarrow\mathbb{R}^{n}$ is Lipschitz continuous and satisfies $\mathbf{f}({\bf 0},{\bf 0}) = {\bf 0}$ and $\nabla_{\bf x}\mathbf{f}({\bf 0},{\bf 0}) = {\bf 0}$.   
%

We define the {\em value function} $v({\bf x}_0) = J({\bf x}^*(\, \cdot\, ;{\bf x}_0),{\bf u}^*(\cdot))$ to be the value of (\ref{eq:oc}) when the optimal control ${\bf u}^*$ and corresponding state ${\bf x}^*$ are found from the initial point ${\bf x}_0$.  Assume that the optimal control is given by the feedback relation
\begin{equation}
\label{eq:feedback}
  {\bf u}(t) = {\mathcal K}({\bf x}(t)).
\end{equation}
For ${\bf f}$, $\ell$, and $v$ smooth enough, and $v$ convex, the feedback relation (\ref{eq:feedback})
satisfies the HJB partial differential equations
\begin{align}
\label{eq:HJB1}
  0 &= \frac{\partial v}{\partial {\bf x}}({\bf x})\left({\bf A}{\bf x} + {\bf B} {\cal K}({\bf x}) +  {\bf f}({\bf x},{\cal K}({\bf x})) \right)
       + \ell({\bf x},{\cal K}({\bf x})), \\
\label{eq:HJB2}
  0 &= \frac{\partial v}{\partial {\bf x}}({\bf x}) \left({\bf B} + \frac{\partial {\bf f}}{\partial {\bf u}}({\bf x},{\cal K}({\bf x})) \right)
       + \frac{\partial \ell}{\partial {\bf u}}({\bf x},{\cal K}({\bf x})).
\end{align}
Ideally, we could solve the HJB equations simultaneously for $v$ and ${\mathcal K}$.  This would provide the desired feedback relation
${\bf u}(t) = {\mathcal K}({\bf x}(t))$.  The value function $v$ can often serve as a Lyapunov function to examine the region of attraction for the controlled system.  Unfortunately, the HJB equations suffer from the curse of dimensionality as they are partial differential equations in $\mathbb{R}^n$.  There have been studies that use \emph{model reduction} to replace the nonlinear dynamics in \eqref{eq:full} with high-fidelity, yet much lower order, reduced dynamics for simple problems~\cite{kunisch2004hjb}.  However, direct solution of the HJB equations will still be a computational challenge even when we can replace the Navier-Stokes equations with modest reduced-order models on the order of 30--50.

Therefore, even with optimal reduced models, solving the nonlinear optimal control problem \eqref{eq:oc}-\eqref{eq:full} is intractable for general nonlinearities. Fortunately, as we outline below, for problems where ${\bf f}$ and $\ell$ have quadratic nonlinearities, there is enough structure to compute nonlinear feedback laws for modest sized problems (at least $n=50$).  This provides a useful tool for reduced systems of flow equations.  Many linear feedback laws are computed from linearized and reduced models of this size in the literature, e.g.~\cite{borggaard2010lfc}.

Our simplification comes from using Kronecker products (which have a long history in the control literature~\cite{brewer1978kronecker,simoncini2016computational}.  Let ${\bf X}\in \mathbb{R}^{i_x\times j_x}$ and ${\bf Y}\in \mathbb{R}^{i_y\times j_y}$, with entries $x_{ij}$ and $y_{ij}$, respectively.  Then ${\bf X}\otimes{\bf Y} \in \mathbb{R}^{i_xi_y \times j_xj_y}$ is the block matrix
\begin{displaymath}
   {\bf X} \otimes {\bf Y} \equiv \left[ \begin{array}{cccc}
     x_{11}{\bf Y} & x_{12}{\bf Y} & \cdots & x_{1j_x}{\bf Y} \\
     x_{21}{\bf Y} & x_{22}{\bf Y} & \cdots & x_{2j_x}{\bf Y} \\
     \vdots  &         &        & \vdots    \\
     x_{i_x1}{\bf Y}& x_{i_x2}{\bf Y}& \cdots & x_{i_xj_x}{\bf Y} \end{array} \right].
\end{displaymath}

\section{The Quadratic-Quadratic Regulator}
To simplify the expressions, we use the Kronecker product description of the quadratic-quadratic regulator (QQR) problem.  Thus, we seek
the control ${\bf u}(t) = \mathcal{K}({\bf x}(t))$ that is the solution to
\begin{equation}\label{eq:quadratic-cost}
  \min_{\bf u} \int_0^\infty {\bf q}_2'\left( {\bf x}(t)\otimes{\bf x}(t) \right) + {\bf r}_2'\left( {\bf u}(t)\otimes{\bf u}(t) \right)\ dt
\end{equation}
subject to
\begin{equation}\label{eq:quadratic-in-state}
  \dot{\bf x}(t) = {\bf A}{\bf x}(t) + {\bf B}{\bf u}(t) + {\bf N}\left( {\bf x}(t)\otimes{\bf x}(t) \right),
\end{equation}
from any ${\bf x}(0) = {\bf x}_0 \in \mathbb{R}^n$.  In this formulation, the matrices above are time-invariant with dimensions
\begin{displaymath}
  {\bf A}\in \mathbb{R}^{n\times n}, \qquad {\bf B}\in \mathbb{R}^{n\times m}, \qquad {\bf N}\in \mathbb{R}^{n\times n^2},
\end{displaymath}
\begin{displaymath}
  {\bf q}_2\in \mathbb{R}^{n^2\times 1}, \qquad \mbox{and} \qquad {\bf r}_2\in \mathbb{R}^{m^2\times 1}
\end{displaymath}
(and $'$ denotes the transpose).
This is merely for a convenient representation for the algorithm below.  Note that we require the standard control systems properties
required by the linear-quadratic regulator (LQR) problem and these can be readily checked with the identities
${\bf q}_2 = {\rm vec}({\bf Q}_2)$ and ${\bf r}_2 = {\rm vec}({\bf R}_2)$ with the usual quadratic cost integrand being ${\bf x}'{\bf Q}_2{\bf x} + {\bf u}'{\bf R}_2{\bf u}$.


We now expand the {\em value function} as
\begin{displaymath}
  v({\bf x}) = \underbrace{{\bf v}_2' \left( {\bf x}\otimes{\bf x} \right)}_{v^{[2]}({\bf x})} + \underbrace{{\bf v}_3' \left( {\bf x}\otimes{\bf x}\otimes{\bf x} \right) }_{v^{[3]}({\bf x})} + \cdots
\end{displaymath}
and the feedback operator as
\begin{displaymath}
  \mathcal{K}({\bf x}) = \underbrace{{\bf k}_1' {\bf x}}_{{\bf k}^{[1]}({\bf x})} + \underbrace{{\bf k}_2' \left( {\bf x}\otimes{\bf x} \right)}_{{\bf k}^{[2]}({\bf x})} + \underbrace{{\bf k}_3' \left( {\bf x}\otimes{\bf x}\otimes{\bf x} \right)}_{{\bf k}^{[3]}({\bf x})} +\cdots.
\end{displaymath}
Note that ${\bf v}_d\in \mathbb{R}^{n^d\times 1}$ and ${\bf k}_d \in \mathbb{R}^{n^d\times m}$. The Hamiltonian-Jacobi-Bellman equations for this problem now has the form
\begin{align}
  \frac{\partial v}{\partial {\bf x}}({\bf x}) \left[ {\bf A}{\bf x} + {\bf B}{\bf u} + {\bf N}({\bf x}\otimes{\bf x})\right] \hspace{1in}\nonumber\\
  \label{eq:HJB1_K}
  +{\bf q}_2'({\bf x}\otimes{\bf x})
  +{\bf r}_2'(\mathcal{K}({\bf x})\otimes \mathcal{K}({\bf x})) &= {\bf 0},\\
  \label{eq:HJB2_K}
  \frac{\partial v}{\partial {\bf x}}({\bf x}) \left[{\bf B}\right] + {\bf r}_2'\mathcal{K}({\bf x}) &= {\bf 0}.
\end{align}

Substituting in the expansions for the value function $v$ and the feedback operator $\mathcal{K}$ into (\ref{eq:HJB1_K}), then collecting 
$O({\bf x}^2)$ terms, we have
\begin{align}\nonumber
  {\bf v}_2' \left( ({\bf Ax}+{\bf B}{\bf k}_1{\bf x})\otimes{\bf x}\! +\! {\bf x}\otimes({\bf Ax}+{\bf B}{\bf k}_1{\bf x}) \right) \\
  + {\bf q}_2'({\bf x}\otimes{\bf x}) + {\bf r}_2'(({\bf k}_1{\bf x})\otimes({\bf k}_1{\bf x})) = {\bf 0}
\end{align}
which, using ${\bf k}_1 = -{\bf R}_2^{-1}{\bf B}'{\bf V}_2$ is equivalent to the algebraic Riccati equation (ARE)
for finding ${\bf V}_2$
\begin{displaymath}
  {\bf A}' {\bf V}_2 + {\bf V}_2{\bf A} - {\bf V}_2 {\bf B}{\bf R}_2^{-1}{\bf B}' {\bf V}_2 + {\bf Q}_2 = {\bf 0}.
\end{displaymath}
Note that it is natural to use the efficient algorithms for solving the ARE and set ${\bf v}_2 = {\rm vec}({\bf V}_2)$.

\subsection{Coefficients of ${\bf v}_{d+1}$}
When we gather higher degree terms, we ignore those terms in (\ref{eq:HJB1_K}) that involve components of $\mathcal{K}$ that have yet to be computed.  Those terms will be updated from (\ref{eq:HJB2_K}).  The degree three terms in (\ref{eq:HJB1_K}) can then be written using the definition ${\bf A}_c = {\bf A}+{\bf Bk}_1$ as
\begin{equation}\label{eq:d2}
\begin{split}
  \left({\bf A}_c \otimes {\bf I}_n \otimes {\bf I}_n + {\bf I}_n \otimes {\bf A}_c \otimes {\bf I}_n  +
  {\bf I}_n \otimes {\bf I}_n \otimes {\bf A}_c\right)' {\bf v}_3 \\
  = -\left({\bf N}'\otimes {\bf I}_n + {\bf I}_n\otimes {\bf N}'\right){\bf v}_2.
\end{split}
\end{equation}
Note that this is a simplification that is independent of ${\bf k}_2$ since collecting those terms then factoring leaves us
with identities following from (\ref{eq:HJB2_K}) involving the ${\bf v}_2$, ${\bf r}_2$ and ${\bf k}_1$ terms that define
${\bf k}_1$.  For example, 
\begin{displaymath}
  (({\bf B}{\bf k}_2)'\otimes {\bf I}_n){\bf v}_2 + ({\bf k}_2'\otimes{\bf k}_1'){\bf r}_2 = 0.
\end{displaymath}
This identity appears in all subsequent collections of similar degree terms since the ${\bf k}_1$ term will always be matched
up with ${\bf k}_d$ terms.  This fact serves to decouple equations
for ${\bf v}_{d+1}$ in (\ref{eq:HJB1_K}) from the equations for ${\bf k}_d$ in (\ref{eq:HJB2_K}).

To write the equations from matching higher degree terms in a more compact way, we define the {\em N-way Lyapunov matrix}
or a special {\em Kronecker sum}~\cite{benzi2015decay} matrix,
\begin{equation}\label{eq:kroneckerSum}
  \mathcal{L}_d({\bf X}) \equiv \underbrace{{\bf X} \otimes \cdots \otimes {\bf I}_n}_{d\ {\rm terms}} + \cdots +\underbrace{{\bf I}_n \otimes \cdots \otimes {\bf X}}_{d\ {\rm terms}}.
\end{equation} 
Then the calculation of ${\bf v}_3$ follows from solving an equation of the form
\begin{equation}
\label{eq:degree3}
  \mathcal{L}_3({\bf A}_c') {\bf v}_3 = -\mathcal{L}_2({\bf N}') {\bf v}_2.
\end{equation}
Once we have ${\bf v}_3$, we can readily compute ${\bf k}_2$ as shown in Section~\ref{sec:k}.  The other terms in the series expansion of the value function lead to equations that have a similar form.
All of the left-hand-sides are generically the same $\mathcal{L}_{d+1}({\bf A}_c') {\bf v}_{d+1}$.  However, 
the right-hand-sides of the equations gather more terms due to the ${\bf r}_2$ term in (\ref{eq:HJB1_K})
and the interactions of the previously computed nonlinear feedback terms with previously computed terms of the value function (that are known and moved to the right-hand-side).  This process is clarified from explicitly writing the next two terms for $v({\bf x})$ below.  For $O({\bf x}^4)$, we have
\begin{equation}
\label{eq:degree4}
  \mathcal{L}_4({\bf A}_c') {\bf v}_4 = -\mathcal{L}_3(({\bf B}{\bf k}_2+{\bf N})') {\bf v}_3 - ({\bf k}_2^\prime \otimes
  {\bf k}_2^\prime ) {\bf r}_2,
\end{equation}
which can be solved for ${\bf v}_4$ once ${\bf k}_2$ is computed from the
solution ${\bf v}_3$ from (\ref{eq:degree3}), 
and
\begin{equation}
\label{eq:degree5}
\begin{split}
  \mathcal{L}_5({\bf A}_c') {\bf v}_5 =& -\mathcal{L}_4(({\bf B}{\bf k}_2+{\bf N})') {\bf v}_4 
  -\mathcal{L}_3(({\bf B}{\bf k}_3)'){\bf v}_3 \\
  &- ( {\bf k}_2\otimes {\bf k}_3 + {\bf k}_3\otimes {\bf k}_2 )' {\bf r}_2.
\end{split}
\end{equation}
Again, once we compute ${\bf k}_3$ from ${\bf v}_4$, we have everything we need to compute ${\bf v}_5$.

In general, while calculation of the coefficients ${\bf v}_{d}$ are described by large linear systems ($\mathcal{L}_{d}({\bf A}_c')\in\mathbb{R}^{n^{d}\times n^{d}}$), there is a great deal of structure.  For example, consistent with the remarks in
\cite{krener2014series}, the eigenvalues of $\mathcal{L}_{d}({\bf A}_c')$ are merely sums of combinations of the eigenvalues
of ${\bf A}_c$.  Therefore, since ${\bf A}_c$ is a stable matrix, $\mathcal{L}_{d}({\bf A}_c')$ will also be.  It is also
immediately obvious that {\em without} the nonlinear term ${\bf N}$ in our state equation, the right-hand-side in (\ref{eq:degree3}) would vanish leading to ${\bf v}_3={\bf 0}$. The remaining equations for ${\bf v}_{d+1}$ would have homogeneous right-hand-sides and thus ${\bf v}_{d+1}={\bf 0}$ for $d=2$ and higher.  This is consistent with the LQR theory.

\subsection{Coefficients of ${\bf k}_{d}$\label{sec:k}}
We now turn our attention to using (\ref{eq:HJB2_K}) to calculate ${\bf k}_d$ from ${\bf v}_{d+1}$.  This is
again straight-forward using the specialized Kronecker sum operator,
\begin{equation}
\label{eq:k_d}
  {\bf k}_d = -{\bf R}_2^{-1} \left(\mathcal{L}_{d+1}({\bf B}') {\bf v}_{d+1}\right)'.
\end{equation}

\subsection{Computing Right-Hand-Side Vectors}
The assembly and solution of linear systems with the form 
\begin{equation}
\label{eq:generic}
  \mathcal{L}_{d+1}({\bf A}_c){\bf v}_{d+1} = {\bf c}
\end{equation}
is  only feasible for small values of $d$ and $n$.  We denote these computations by {\em full Kronecker} in 
Section~\ref{sec:numerical}.  The advantage of the Kronecker product structure is that we can perform operations with Kronecker product matrices {\em without} actually forming the large block matrix.  The main issue that we deal with in this section is calculating the terms on the right-hand-sides of e.g. (\ref{eq:degree3})--(\ref{eq:degree5}) or (\ref{eq:k_d}). Solution of the system (\ref{eq:generic}) is described in the next section. 

To calculate ${\bf c}$ for (\ref{eq:degree3})--(\ref{eq:k_d}) involves two types of terms.  The first involves the multiplication of a Kronecker form with a vector ${\bf r}_2$.  Recall, e.g.~\cite{brewer1978kronecker}, that 
\begin{equation}\label{eq:kron*v}
  ({\bf X}\otimes{\bf Y}){\bf r}_{2} = {\rm vec}({\bf Y}'{\bf R}{\bf X}),
\end{equation}
where ${\bf R}$ has the appropriate dimensions and ${\bf r}_{2} = {\bf vec}({\bf R})$.  Therefore, the terms involving ${\bf r}_2$ only require
matrix multiplications and no assembly of the Kronecker product is required.

The second type of term are products of the Kronecker sum with a ${\bf v}_{d+1}$: $\mathcal{L}_{d+1}({\bf X}){\bf v}_{d+1}$.
Using the definition of (\ref{eq:kroneckerSum}), we have to calculate $d+1$ different multiplications 
of the Kronecker products with ${\bf v}_{d+1}$.  Using the fact
that the Kronecker product is associative, we write
\begin{displaymath}
  {\bf I}_{n^\ell} = \underbrace{{\bf I}_n\otimes \cdots\otimes {\bf I}_n}_{\ell {\rm terms}}.
\end{displaymath}
The multiplications can be reduced to three different cases
\begin{displaymath}
  ({\bf X} \otimes {\bf I}_{n^d}){\bf v}_{d+1}, \ \ ({\bf I}_{n^{d-\ell}}\otimes {\bf X}\otimes {\bf I}_{n^{\ell}}){\bf v}_{d+1},
  \ \ \mbox{and} \ \ ({\bf I}_{n^d}\otimes {\bf X}){\bf v}_{d+1}.
\end{displaymath}
Here the relation (\ref{eq:kron*v}) and the associative law for Kronecker products are useful.  The first and last terms
above can be handled by the appropriate reshaping of ${\bf v}_{d+1}$ and multiplying with ${\bf X}$ (the multiplication by ${\bf I}_{n^d}$ is trivial).  The associative law allows us to handle all of the intermediate terms recursively as
\begin{align*}
   ({\bf I}_{n^{d-\ell}}\otimes {\bf X}\otimes {\bf I}_{n^{\ell}}){\bf v}_{d+1} &= (({\bf I}_{n^{d-\ell}}\otimes {\bf X})\otimes {\bf I}_{n^{\ell}}){\bf v}_{d+1}\\
       &= ({\bf I}_{n^{d-\ell}}\otimes ({\bf X}\otimes {\bf I}_{n^{\ell}})){\bf v}_{d+1}.
\end{align*}
The grouping can be done to maximize the size of the free identity matrix.

\subsection{Linear System Solutions}
The Kronecker structure leads to larger systems (\ref{eq:generic}), but are now ameneble to modern high performance algorithms~\cite{kolda2009tensor,chen2019recursive,simoncini2016computational}.  Many of these algorithms,
e.g. \cite{chen2019recursive} utilize a real Schur factorization of the matrix ${\bf A}_c$.
For this study, we used the {\em recursive} algorithms in \cite{chen2019recursive} for Laplace-like equations. Their software was 
trivially modified to take advantage of the fact that the same term ${\bf A}_c$ appears in every block and gave the
system exactly the form (\ref{eq:kroneckerSum}).

\section{Numerical Results\label{sec:numerical}}
We present two sets of results.  The first is a challenging verification test using randomly generated matrices.  This
is challenging since the conditioning of these systems are poor for many of our randomly generated samples.  Our study
is reproducible since we reset the random seed before each test.  The second set a discretized control problem involving
the one-dimensional Burgers equation.  The systems are much better conditioned in this case and there is a notion of
convergence as the problem sizes in our tests grow.

\subsection{Random System Study}
In this section, we perform computational tests for both performance and accuracy.  For 
accuracy, we compare against the feedback matrices computed using the Nonlinear Systems Toolbox (NST)~\cite{krener2015NST}.
We note that the NST is designed for a wide range of control problem formulations and general nonlinearities.  Thus it is not used as a performance measure but rather to provide a sense of general speed and accuracy.  Since much of the use of the symbolic toolbox in NST is in the preprocessing step, we remove this calculation from our comparative timings (though provide those times separately for completeness).

For reproducibility, we provide the source to generate our numerical findings in Matlab below.
\begin{verbatim}
  rng(0,'v5uniform'); % set random seed
  A = rand(n,n);
  B = rand(n,m);
  N = rand(n,n*n);
  Q = eye(n); 
  R = eye(m);
\end{verbatim}
All computations were performed on a 2017 Macbook Pro with a 3.1GHz Inter Core i7 processor and 16GB of RAM
using MATLAB version R2019b.

Our first set of tests computed the degree 2 feedback term, ${\bf k}_2$, and the required ${\bf v}_3$ component of
the value function.  Problem sizes from 6 to 20 were randomly generated and the Matlab runtimes are
reported in Table~\ref{table:d2Times}.  Generally speaking, only the first significant digits or two of the CPU times 
are meaningful since we only averaged the time over a small set of experimental runs.  The trend is clear, however, 
that our QQR algorithm is exceptionally fast primarily due to the careful work assembling the right-hand-sides and the
recursive blocked algorithms of Chen and Kressner~\cite{chen2019recursive}.  
The computational trend for the QQR is even more pronounced when we
extend these tests to include both the degree 3 feedback term, ${\bf k}_3$, and the corresponding ${\bf v}_4$ 
component of the value function (seen in Table~\ref{table:d3Times}).  The calculation of the full Kronecker problem for $n=16$ encountered a 
memory limit error, so computations beyond that could not be carried out.  However, the recursive solver and
specially tailored matrix-vector products allowed us to continue calculations well beyond the limits seen
with the full Kronecker form.  In fact, we were able to solve an order $n=40$ random system
with a degree $d=3$ feedback law in 12.72 seconds and an order $n=30$ random system with a
degree $d=4$ feedback law in 142.26 seconds.  This demonstrates that modest size problems in the QQR
framework can be readily incorporated into control design workflows.
\begin{table}[h]
\caption{Random: CPU Time for Degree 2 Terms\label{table:d2Times}}
\begin{center}
\begin{tabular}{|c||c|c|cr@{.}l|}
\hline
n & recursive & full Kronecker & \multicolumn{3}{|c|}{NST \ \ \ (symbolic calc.)} \\
\hline
 6 & 0.03584  & 0.00152  & 0.0123 &(0&5665) \\
 8 & 0.03958  & 0.00552  & 0.0258 &(0&9183) \\
10 & 0.03531  & 0.02020  & 0.0515 &(1&4596) \\
12 & 0.04466  & 0.07945  & 0.1006 &(2&2126) \\
14 & 0.05852  & 0.22527  & 0.2528 &(3&6933) \\
16 & 0.06112  & 0.67803  & 0.5932 &(5&6935) \\
18 & 0.07812  & 1.55314  & 1.2353 &(8&8498) \\
20 & 0.09467  & 3.50995  & 2.5385 &(13&3123) \\
\hline
\end{tabular}
\end{center}
\end{table}
\begin{table}[h]
\caption{Random: Cumulative CPU Time for Degree 3 Terms\label{table:d3Times}}
\begin{center}
\begin{tabular}{|c||c|r@{.}l|r@{.}lr@{.}l|}
\hline
n & recursive & \multicolumn{2}{|c|}{full Kronecker} & \multicolumn{4}{|c|}{NST \ \ \ (symbolic calc.)} \\
\hline
 6 & 0.00700  & 0&04937  & 0&0333 &(1&6310) \\
 8 & 0.04363  & 0&87317  & 0&1219 &(4&0723) \\
10 & 0.06942  & 6&6439   & 0&4928 &(9&1590) \\
12 & 0.09098  &53&2562   & 1&9644 &(20&2047) \\
14 & 0.21810  &588&826   & 7&1282 &(41&0284) \\
16 & 0.53792  &\multicolumn{2}{|c|}{\rm not computed} & 23&9001 &(85&0918) \\
18 & 0.63851   &\multicolumn{2}{|c|}{\rm not computed} & 68&1008 &(174&518) \\
20 & 0.82421   &\multicolumn{2}{|c|}{\rm not computed} &170&217  &(352&068) \\
\hline
\end{tabular}
\end{center}
\end{table}

There is no truth model for testing the accuracy of the HJB series solution, so we relied
on the well-tested NST software to compare against.  Therefore, relative ${\ell^2}$ 
errors are those reported with respect to NST solutions.  This was a 
challenge since NST computes its solutions in a compact Taylor series format 
(returning coefficients of the unique monomial terms) whereas the Kronecker 
formulation introduces a lot of repeated monomials.  Thus there are multiple correct
representations for the coefficients (even though only one unique representation
will be calculated).  To compare to the coefficients in 
compact Taylor series form required us to accumulate all of the coefficients for equivalent
monomials to generate the comparisons.  The relative errors in $\ell^2$ 
are those of the summed coefficients (i.e. summing the coefficients for $x_1x_2$ with those
from $x_2x_1$ etc.).  Tables~\ref{table:d2errors} and
\ref{table:d3errors} show very good agreement for small problems.
Note that both the full Kronecker system and the NST systems generated 
condition number warnings for $n=16$ and $n=20$ for the $d=2$ solutions and
the {\em NST} generated an additional warning when $d=3$ for the $n=18$ case.  
These warnings did not occur in the recursive blocked algorithms since the full system was never assembled.
However, the systems that are solved in the recursive and full Kronecker columns
are the same, so many of the large relative errors can be explained by ill-conditioning. 
Therefore, our next numerical example in Section~\ref{sec:burgers} will have more desirable control-theoretic
properties.
\begin{table}[h]
\caption{Random: Relative Errors in $k^{[2]}$ and $v^{[3]}$\label{table:d2errors}}
\begin{center}
\begin{tabular}{|c||cc|cc|}
\hline
  & \multicolumn{2}{|c|}{recursive} & \multicolumn{2}{|c|}{full Kronecker} \\
n & error $k^{[2]}$ & error $v^{[3]}$ & error $k^{[2]}$ & error $v^{[3]}$ \\
\hline
 6 & 4.09e-12 & 3.09e-12  & 3.88e-12 & 2.61e-12   \\
 8 & 1.70e-11 & 1.67e-11  & 2.05e-11 & 2.66e-11   \\
10 & 4.92e-08 & 4.95e-08  & 1.25e-06 & 1.27e-06   \\
12 & 1.19e-10 & 1.04e-10  & 3.93e-10 & 3.72e-10   \\
14 & 9.26e-07 & 9.95e-07  & 7.42e-06 & 7.88e-06   \\
16 & 1.06e-04 & 1.13e-04  & 3.80e-03 & 3.97e-03   \\
18 & 1.44e-04 & 1.50e-04  & 5.48e-03 & 5.75e-03   \\
20 & 2.65e-02 & 3.23e-02  & 1.68e+00 & 1.70e+00   \\
\hline
\end{tabular}
\end{center}
\end{table}
\begin{table}[h]
\caption{Random: Relative Errors in $k^{[3]}$ and $v^{[4]}$\label{table:d3errors}}
\begin{center}
\begin{tabular}{|c||cc|cc|}
\hline
  & \multicolumn{2}{|c|}{recursive} & \multicolumn{2}{|c|}{full Kronecker} \\
n & error $k^{[3]}$ & error $v^{[4]}$ & error $k^{[3]}$ & error $v^{[4]}$ \\
\hline
 6 & 1.68e-10 & 1.11e-10 & 3.37e-11 & 2.71e-11 \\
 8 & 9.18e-10 & 6.85e-10 & 1.17e-09 & 1.00e-09 \\
10 & 7.52e-04 & 9.11e-04 & 2.32e-03 & 2.81e-03 \\
12 & 7.21e-10 & 1.04e-09 & 1.15e-07 & 1.28e-07 \\
14 & 7.95e-04 & 5.80e-04 & 1.86e-02 & 1.38e-02 \\
16 & 2.11e-01 & 2.61e-01 & {\rm not computed}  & {\rm not computed} \\
18 & 7.90e+00 & 1.15e+00 & {\rm not computed}  & {\rm not computed} \\
\hline
\end{tabular}
\end{center}
\end{table}

\subsection{Burgers Equation\label{sec:burgers}}
As a more structured test problem, we consider the QQR problem with a discretization
of the Burgers equation.  This test problem has a long history in the study of control for 
distributed parameter systems, e.g.~\cite{thevenet2009nonlinear}, including the 
development of effective computational methods, e.g.~\cite{burns1990control}.  Thus, we
consider 

We consider the specific problem found in \cite{borggaard2018computation} but with
two control inputs ($m=2$) that consist of uniformly distributed sources over 
disjoint patches.  Thus, we have a bounded input operator.  The formal description of
the problem is 
\begin{displaymath}
  \min_{\bf u}J(z,u) = \int_0^\infty \left(\int_0^1 z^2(\xi,t)\ d\xi + {\bf u}'(t){\bf u}(t)\right)\ dt
\end{displaymath}
subject to
\begin{align*}
  \dot{z}(x,t) &= \epsilon z_{xx}(x,t) - \frac{1}{2}\left(z^2(x,t)\right)_x 
  \\
  &\hspace{.4in}+ \sum_{k=1}^m 
  \chi_{[(k-1)/m,\ k/m]}(x) u_k(t)\\
  z(\cdot,0) &= z_0(\cdot) \in H_{\rm per}^1(0,1),
\end{align*}
where $\chi_{[a,b]}(x)$ is the characteristic function over $[a,b]$. We discretized the state equations with $n$ linear finite elements, set $m=2$, and chose $\epsilon=0.001$ to make the nonlinearity significant.  

The discretized system fits within the QQR framework (\ref{eq:quadratic-cost})-(\ref{eq:quadratic-in-state}).  The matrices ${\bf A}$, ${\bf B}$ and ${\bf N}$ come from the finite element approximation.  The matrix ${\bf Q}_2$ is the finite element mass matrix and the matrix ${\bf R}_2={\bf I}_m$.


With the same rationale, the NST solution is used as the truth model and used to
compute relative errors.  It is evident that a well-conditioned feedback control problem
produces much more consistent errors.  As expected, the $n=16$ case also breaks down here
(more memory is required since $m=2$), but the accuracy of the recursive block algorithm
is easily observed in Tables~\ref{table:bd2errors} and \ref{table:bd3errors}.  There is
very little growth in the discrepancy between the two solutions with increasing values of $n$.
This suggests a convergence to an infinite-dimensional representation that we will investigate
in future studies.
\begin{table}[h]
\caption{Burgers: Relative Errors in $k^{[2]}$ and $v^{[3]}$\label{table:bd2errors}}
\begin{center}
\begin{tabular}{|c||cc|cc|}
\hline
  & \multicolumn{2}{|c|}{recursive} & \multicolumn{2}{|c|}{full Kronecker} \\
n & error $k^{[2]}$ & error $v^{[3]}$ & error $k^{[2]}$ & error $v^{[3]}$ \\
\hline
10 & 6.34e-12 & 3.99e-12  & 6.34e-12 & 3.99e-12   \\
12 & 1.02e-11 & 4.94e-12  & 1.02e-11 & 4.94e-11   \\
14 & 2.53e-11 & 1.85e-11  & 2.53e-11 & 1.85e-11   \\
16 & 2.33e-11 & 1.24e-11  & {\rm not computed} & {\rm not computed} \\
18 & 3.63e-11 & 1.76e-11  & {\rm not computed} & {\rm not computed} \\
20 & 7.43e-11 & 4.44e-11  & {\rm not computed} & {\rm not computed} \\
\hline
\end{tabular}
\end{center}
\end{table}
\begin{table}[h]
\caption{Burgers: Relative Errors in $k^{[3]}$ and $v^{[4]}$\label{table:bd3errors}}
\begin{center}
\begin{tabular}{|c||cc|cc|}
\hline
  & \multicolumn{2}{|c|}{recursive} & \multicolumn{2}{|c|}{full Kronecker} \\
n & error $k^{[3]}$ & error $v^{[4]}$ & error $k^{[3]}$ & error $v^{[4]}$ \\
\hline
10 & 1.21e-12 & 3.02e-12 & 1.21e-12 & 3.02e-12 \\
12 & 1.74e-11 & 4.79e-12 & 1.74e-11 & 4.80e-11 \\
14 & 3.66e-11 & 1.67e-11 & 3.66e-11 & 1.67e-11 \\
16 & 3.07e-11 & 1.00e-11 & {\rm not computed} & {\rm not computed} \\
18 & 4.76e-11 & 1.52e-11 & {\rm not computed} & {\rm not computed} \\
20 & 1.22e-10 & 4.60e-11 & {\rm not computed} & {\rm not computed} \\
\hline
\end{tabular}
\end{center}
\end{table}

Instead of presenting both degree 2 and the accumulation up to degree three,
we only report the accumulated times up to degree 3 in Table~\ref{table:bd3Times}.
Again, the computational speed to solve the QQR problem for the degree 2 and
degree 3 feedback terms is impressive.
\begin{table}[h]
\caption{Burgers: Cumulative CPU Time for Degree 3 Terms\label{table:bd3Times}}
\begin{center}
\begin{tabular}{|c||c|r@{.}l|r@{.}lr@{.}l|}
\hline
n & recursive & \multicolumn{2}{|c|}{full Kronecker} & \multicolumn{4}{|c|}{NST \ \ \ (symbolic calc.)} \\
\hline
10 & 0.07355  & 6&2939   & 0&5367 &(12&0943) \\
12 & 0.06801  &49&9701   & 2&2008 &(26&7163) \\
14 & 0.19543  &474&602   & 7&4379 &(50&843) \\
16 & 0.39385  &\multicolumn{2}{|c|}{not computed} & 25&9467 &(108&846) \\
18 & 0.42039  &\multicolumn{2}{|c|}{not computed} & 69&2946 &(203&07) \\
20 & 0.50770  &\multicolumn{2}{|c|}{not computed} &179&731  &(388&069) \\
\hline
\end{tabular}
\end{center}
\end{table}

To test the possibility of using this software for larger problems, we increased the degree and order
to see what sizes could be computed in two to three minutes for this problem.  We found that we could
solve a degree $d=4$ feedback law for order $n=32$ and $m=2$ in 182.9 seconds.
Furthermore, a degree $d=3$ feedback law for order $n=64$ and $m=2$ was computed in 130.2 seconds.
Thus, some higher degree control laws are feasible with modest discretizations in one-dimensional
problems.  It is reasonable to expect that this would be an attractive option to try when developing
feedback control laws from modest reduced models of fluid systems.

\section{Conclusions and Future Work}
We presented a special formulation of the Al'Brekht polynomial approximation for the
quadratic-quadratic regulator problem.  Writing the expansions in terms of Kronecker
products leads to a series of progressively larger linear systems for the next terms
in the expansion.  While easy to write down and implement, efficiency is only achieved
by exploiting new numerical linear algebra tools that avoid the assembly of the large,
dense systems~\cite{kolda2009tensor,chen2019recursive}.  Repeatable numerical experiments
with random linear systems of different order confirm the efficiency of our approach.
We performed a comparison with a general, well-developed software tool, the Nonlinear Systems
Toolbox~\cite{krener2015NST}, to verify our implementation.  Our solution method was 
competitive with NST in terms of CPU time even if we neglect the overhead in using Matlab's
symbolic toolbox (we described an effective means to compute the derivatives of the
system that are required by NST using automatic differentiation in a previous 
paper~\cite{borggaard2018computation}).

Our future work will evolve down three paths.  The first will be to better understand the
numerical trade-offs between our use of \cite{chen2019recursive} for our application and
a linear system that is built for a compact Taylor series representation of the problem.
This compact Taylor series removes redundant variables (e.g. coefficients of $x_1 x_2$ and 
$x_2 x_1$ can be combined).  Since NST uses the compact Taylor series approach, we have
already built the restriction and prolongation operators between polynomials up to order
5 for our verification.  Part of this will include a better understanding of the discrepancies
between the solutions for some of our randomly generated test cases.  

Our second path is to consider natural generalizations to include in this software framework.
This would include the investigation of more general control costs, addition of descriptor 
systems~\cite{xu1993hamilton}, and the related observer problem.

Finally, we will apply this to more significant applications than the one-dimensional Burgers
equation.  In particular, study how this work could be used in conjunction with reduced
models of complex flows that result in quadratic-in-state systems.

This software is available for download at {\tt https://github.com/jborggaard/QQR.git}.

\bibliographystyle{plain}
\bibliography{refs}

\end{document}